\input amssym.def
\input amssym.tex
\hsize=140mm \vsize=270mm \magnification=1030 
\input xy
\xyoption{all}

\centerline{\bf ORBIFOLDES, VARI\'ET\'ES SP\'ECIALES ET CLASSIFICATION}

\centerline{ \bf  DES VARI\'ET\'ES 
K\" AHL\'ERIENNES COMPACTES}

\

\centerline{\bf F.Campana}

\

Le pr\'esent texte est une introduction \`a [C01] , 
o\`u l'on trouvera les d\'etails techniques. On suppose le lecteur familier avec les 
notions de base de la g\'em\'etrie complexe (applications m\'eromorphes, formes diff\'erentielles 
holomorphes, fibr\'es en droites amples, courbes projectives complexes).

\

{\bf 0. INTRODUCTION.}

\

{\bf Notations:} On notera $X$ une vari\'et\'e complexe compacte et connexe K\"ahl\'erienne 
de dimension complexe $n$, 
$K_X:=det(\Omega_X^1) $ son fibr\'e canonique, dont les sections locales sont les formes volumes holomorphes.

Une {\bf fibration} sera une application m\'eromorphe surjective \`a fibres connexes $f:X\rightarrow Y$, 
avec $Y$ un espace analytique complexe normal de dimension $p\leq n$. On notera $X_y$ sa fibre ``g\'en\'erale" (ie: au-dessus de
 $y\in Y$, dans une intersection d\'enombrable d'ouverts de Zariski denses. Les fibres 
de $f$ sont les images de celles obtenues 
apr\`es r\'esolution des ind\'eterminations). Deux fibrations sont \'equivalentes si elles ont la m\^eme famille de fibres, apr\`es identification bim\'eromorphe des vari\'et\'es $X,X'$ sur lesquelles elles sont d\'efinies.

\

{\bf PROBL\`EME:} ``Classifier" $X$, \`a \'equivalence bim\'eromorphe pr\`es.

\

Trois g\'eom\'etries ``pures" (``sph\'erique", ``plate" et ``hyperbolique") 
se d\'egagent naturellement, d\'efinies par le ``signe" (n\'egatif, nul, ou positif) de $K_X$, 
ou de mani\`ere \'equivalente, par le signe de la courbure de Ricci de $X$, avec un signe oppos\'e.

\

Dans le cas bien connu des courbes ($n=1$), ces trois g\'eom\'etries sont d\'etermin\'ees par un invariant 
topologique: le genre $g(X)\geq 0$. On a la trichotomie suivante:

\

{\bf g\'eom\'etrie ``sph\'erique"} si $K_X<0$, ie: ssi $Ricci(X)>0$, ssi $g(X)=0$, 
ssi $X\cong \Bbb P^1(\Bbb C)$.

\

{\bf g\'eom\'etrie (Ricci-)``plate"} si $K_X=0$, ie: ssi $Ricci(X)=0$, ssi 
$g(X)=1$, ssi $X\cong \Bbb C/ \Lambda,\Lambda$ un r\'eseau de $\Bbb C$, ssi $X$ 
est une courbe elliptique. 

\

{\bf g\'eom\'etrie ``hyperbolique"} si $K_X>0$, ie: ssi $Ricci(X)<0$, ssi 
$g(X)\geq 2$, ssi $X\cong \Bbb D/ \Gamma,\Gamma$ un r\'eseau cocompact de $SU(1,1)$, et 
$\Bbb D$ le disque unit\'e de $\Bbb C$.

\

Ces trois g\'eom\'etries sont distingu\'ees aussi par d'autres 
invariants usuels (groupe fondamental, pseudom\'etrique de Kobayashi, 
r\'epartition des points $K$-rationnels,...). Les deux premi\`eres g\'eom\'etries ont des comportements 
qualitativement voisins, antith\'etiques de celui de la troisi\`eme. Par exemple, la 
pseudom\'etrique de Kobayashi est nulle pour les g\'eom\'etries plate et sph\'erique, mais une m\'etrique 
dans le cas hyperbolique; de mani\`ere analogue, si $X$ est d\'efinie 
sur un corps de nombres $K$ assez grand, l'ensemble $X(K)$ de ses points $K$-rationnels est dense dans les cas 
plat et sph\'erique, mais est fini dans le cas hyperbolique, par le th\'eor\`eme de Faltings, conjectur\'e par 
Mordell.

On donnera ci-dessous {\bf (\S 2.C)}
des g\'en\'eralisations bim\'eromorphes naturelles de ces g\'eom\'etries en toute dimension. Le probl\`eme de 
classification comporte donc une premi\`ere \'etape: 

la classification (et l '\'etude) des vari\'et\'es $X$ 
de g\'eom\'etrie ``pure". Cette \'etape ne sera pas abord\'ee ici. 
L'\'etude des g\'em\'etries pures ne suffit pas \`a r\'esoudre le probl\`eme de classification pos\'e.

En effet, d\`es que $n\geq 2$, les produits de courbes montrent que 
le signe de $K_X$ n'est pas d\'efini en g\'en\'eral, et 
peut \^etre mixte, pr\'esentant des directions n\'egatives, nulles ou positives en chaque point.

\

{\bf OBJECTIF:} D\'ecomposer $X$, compacte K\" ahler 
arbitraire, en ses ``composantes" sph\'erique, plate et 
hyperbolique, de mani\`ere \`a r\'eduire (dans une large mesure) 
la compr\'ehension globale de la structure des $X$ arbitraires \`a celle des $X$ de g\'eom\'etrie pure. 
Une telle 
d\'ecomposition ne peut prendre la forme simple de produits. En g\'en\'eral, on devra se 
contenter de fibrations 
dans lesquelles base et fibres ont une g\'eom\'etrie pure donn\'ee.
Plus pr\'ecis\'ement, on va: 

\

{\bf 1.} D\'efinir {\bf (\S 6.A)} une g\'eom\'etrie {\bf ``sp\'eciale"} par l'absence sur $X$ de 
fibration \`a base hyperbolique au sens {\bf ``orbifolde"}. Les vari\'et\'es plates et sph\'eriques sont sp\'eciales.
L'introduction d'une structure 
orbifolde sur la base d'une fibration est en fait l'ingr\'edient nouveau, n\'ecessaire, mais suffisant, pour 
b\^atir une th\'eorie coh\'erente de la d\'ecomposition des vari\'et\'es $X$. La classification doit 
donc \^etre \'etendue \`a la cat\'egorie des orbifoldes, \`a laquelle la plupart des techniques 
existantes s'\'etendent 
sans difficult\'e.

\

{\bf 2.} Construire (au {\bf (\S 7.A)}), pour toute $X$, une fibration, fonctorielle en $X$,
$c_X:X\rightarrow C(X)$ (le {\bf coeur} de $X$), 
caract\'eris\'ee par le fait que:

a. ses fibres g\'en\'erales sont sp\'eciales.

b. sa base orbifolde est hyperbolique.

Cette fibration d\'ecompose donc $X$ en ses composantes sp\'eciale (les fibres de $c_X$) 
et hyperbolique (la base orbifolde de $c_X$). Et $X$ est sp\'eciale ssi $C(X)$ est r\'eduit \`a un point.

\

{\bf 3.} D\'ecomposer ensuite (au {\bf (\S 7.B)}) $c_X$ fonctoriellement comme une tour de fibrations 
\'el\'ementaires de deux types: $r$ (resp. $J$), 
dont les fibres sont sph\'eriques (resp. plates) dans un sens orbifolde ad\'equat. Le r\'esultat final 
prend la forme: $c_X=(J\circ r)^n$.

Lorsque $X$ est sp\'eciale, cette d\'ecomposition du coeur exprime donc $X$ comme tour de fibrations \`a
 fibres orbifoldes alternativement sph\'eriques et plates. La g\'eom\'etrie sp\'eciale apparait ainsi 
comme la combinaison orbifolde 
des deux premi\`eres g\'eom\'etries, sph\'erique et plate, de la trichotomie.

\

{\bf 4.} Donner {\bf (\S 8)}, pour $X$ compacte K\" ahler, \`a l'aide du ``coeur" 
une d\'ecomposition qualitative 
conjecturale tr\`es simple de la pseudom\'etrique 
de Kobayashi $d_X$ de $X$, et de l'ensemble $X(K)$ des
 points $K$-rationnels de $X$ (si $X$ est projective, 
d\'efinie sur un corps de nombres $K$ assez grand). Ces conjectures \'etendent celles de S.Lang, 
qui traitent du cas particulier o\`u $X$ est ``hyperbolique".

A titre d'exemple, on conjecture que $X$ est sp\'eciale {\bf si et seulement si} $d_X\equiv 0$, 
{\bf si et seulement si} $X(K)$ est dense dans $X$ (pour $X$ comme ci-dessus, et $K$ assez grand). Ces assertions pr\'esentent 
ainsi la g\'eom\'etrie sp\'eciale comme la combinaison des g\'eom\'etries sph\'erique et plate 
\'egalement des points de vue arithm\'etique et Kobayashi pseudom\'etrique.

\

La trichotomie initiale de la g\'eom\'etrie complexe se r\'eduit ainsi 
\`a une dichotomie sp\'eciale-hyperbolique 
de nature plus fondamentale, semble-t-il.

\

{\bf 1. DIMENSION CANONIQUE (OU: ``DE KODAIRA")}

\

{\bf A. Dimension de Kodaira-Iitaka d'un fibr\'e en droites.} 

\

Soit $X$ une vari\'et\'e projective complexe de dimension $n$, et $L$ un fibr\'e holomorphe 
en droites complexes sur $X$ (par exemple: $L=K_X$, si $X$ est lisse). On va d\'efinir une 
notion de ``positivit\'e g\'eom\'etrique" pour $L$ \`a l'aide du comportement asymptotique de 
$h^0(X,mL):=dim_{\Bbb C}(H^0(X,mL))$ quand $m>0$ tend vers $+\infty$. (On utilise la notation 
additive $mL:=L^{\otimes m}$). L'id\'ee est que les $h^0(X,mL)$  croissent 
(\`a $(X,L)$ fix\'e) comme un polyn\^ome dont le degr\'e est la dimension de la base d'une fibration 
dont les fibres sont les sous-vari\'et\'es maximales de $X$ le long desquelles $L$ est ``plat".

\

{\bf 1.1. Exemples:}

{\bf 1.} si $L=\cal O$$_X$ (le fibr\'e trivial), on a: $h^0(X,mL)=1, \forall m>0$.

{\bf 2.} Si $L$ est ample (ie: admet une m\'etrique hermitienne de courbure positive), alors:
 
$h^p(X,mL)=0,\forall p>0,\forall m>m_0$ (annulation de Kodaira). 

Donc:
 $h^0(X,mL)=\chi(X,mL)\sim C(X,L).m^n$, quand ${m\rightarrow +\infty}$, o\`u $C(X,L)=L^n/n!>0$.

{\bf 3.} Par contre, si $-L$ est effectif, c'est-\`a-dire si $L$ a une section m\'eromorphe non-nulle 
$s$ ayant des p\^oles, mais pas de z\'ero, alors $h^0(X,mL)=0$ pour tout $m>0$ (car $(m-m)L=\cal O$$_X$ 
aurait alors une section non nulle ayant des z\'eros d'ordre $m$ l\`a o\`u $s$ a des p\^oles).

{\bf 4.} On peut enfin remarquer que si $f:X\rightarrow Y$ est une 
fibration holomorphe, avec $p:=dim_{\Bbb C}(Y)$, 
et si $M$ est un fibr\'e en droites holomorphe sur $Y$, alors: $h^0(X,mf^*(M))=h^0(Y,mM),\forall m>0$. 

Si $M$ est ample sur $Y$, on a donc: $h^0(X,mf^*(M))\sim C(Y,M)m^p$, quand ${m\rightarrow +\infty}$.

{\bf 5.} Si $h^0(X,mL)\sim A.m^p$, quand ${m\rightarrow +\infty}$, on a aussi: 
$Log ( h^0(X,mL))\sim p.Log ( m)$, d'o\`u la:

\

{\bf 1.2. D\'efinition:} Soit $X$ un espace analytique complexe compact 
irr\'eductible, et $L$ un fibr\'e en droites holomorphe sur $X$. On pose: 
$\kappa(X,L):=limsup_{m\rightarrow +\infty}(Log(h^0(X,mL))/Log(m))$. On appelle 
cet invariant la {\bf dimension de Kodaira-Iitaka} de $L$  

\

{\bf 1.3. Exemples:}  Si $(X,L)$ est comme dans les exemples {\bf 1,2,3,4,5} pr\'ec\'edents,
 $\kappa(X,L)$  vaut donc respectivement: $0,n,-\infty,p$. 

\

Plus g\'en\'eralement, on montre les r\'esultats suivants (voir [U]): 

\

{\bf 1.4. Proposition:} 

{\bf 1.} $\forall m>0: \kappa (X,L)=\kappa (X,mL)\in \lbrace -\infty,0,1,...,n\rbrace$.

{\bf 2.} $\kappa (X,L)=-\infty$ ssi $h^0(X,mL)=0,\forall m>0$.

{\bf 3.} $\kappa (X,L)=d\geq 0$ ssi $Bm^d\geq h^0(X,mL)\geq Am^d,\forall m>m_0$, 
pour deux constantes $B>A>0$, quitte \`a remplacer $L$ par un multiple $m_1L$ ad\'equat.

{\bf 4.} Si $\kappa(X,L)=d\geq 0$, il existe une (\`a \'equivalence bim\'eromorphe pr\`es unique)
 fibration $\Phi_L:X\rightarrow V_L$ telle que $dim(V_L)=d$, et $\kappa(X_v,L_v)=0$, si $X_v$ est la fibre 
g\'en\'erale de $\Phi_L$, et $L_v$ la restriction de $L$ \`a $X_v$. Cette fibration est la {\bf fibration 
d'Iitaka} associ\'ee \`a $L$.

{\bf 5.} Si $f:X\to Y$ est une fibration, et $L$ un fibr\'e en droites sur $X$, alors on a l'in\'egalit\'e 
dite (``additivit\'e facile"): $\kappa(X,L)\leq p+\kappa(X_y,L_{\vert X_y})$, avec $p=dim(Y)$.

\

{\bf 1.5. Remarques:} L'assertion {\bf 4.} pr\'ec\'edente est un substitut (tr\`es faible) du cas
$L=\Phi^*(M)$, avec $M$ ample sur $V_L$, comme dans l'exemple {\bf (1.1.4)}. En effet: 
l'application $\Phi _L$ de l'assertion {\bf 1.4.4} ci-dessus est tr\`es classique: 
elle n'est autre que l'application $\Phi_{mL}:X\to \Bbb P(H^0(X,mL)^*):=\Bbb P^N$ associ\'ee au 
syst\`eme lin\'eaire $\vert mL\vert$, pour $m$ grand et divisible. Cette application associe \`a 
$x\in X$ g\'en\'erique l'hyperplan constitu\'e des sections qui s'annulent en $x$. Si $L$ est ample, 
cette application est un plongement de $X$ dans $\Bbb P^N$, et les sections de $mL$ s'annulent 
sur les intersections de $\Phi_{mL}(X)$ avec les hyperplans de $\Bbb P^N$.

Si $L=f^*(M)$, comme dans l'exemple {\bf (1.1.4)} ci-dessus, alors $\Phi_M$ est 
un plongement de $Y$, et clairement, $\Phi_L=f\circ \Phi_M$, puisque $H^0(X,mL)=f^*(H^0(Y,mM))$. 
Donc $\Phi_L(X)=Y$ est bien de dimension $p$, et 
$\Phi_L=f$.

\

{\bf B. La Dimension Canonique.} 

\

{\bf 1.6. D\'efinition:} Soit $X$ vari\'et\'e analytique complexe compacte connexe (lisse). On pose:

$\kappa(X):=\kappa(X,K_X)\in \lbrace -\infty,0,1,...,n\rbrace$. 

Si $\kappa(X)=n(=dim(X)$, on dit que 
$X$ est de {\bf type g\'en\'eral}.

On appelle $\kappa(X)$ la 
{\bf dimension canonique (ou: ``de Kodaira") de $X$}.

(Le terme de ``dimension canonique est d\^u \`a Moishezon ([Mo], o\`u il l'a introduite). 
Cette notion apparait en 
fait pour la premi\`ere fois dans [Sh], semble-t-il, mais dans le cas particulier des surfaces. 
Le terme, usuel, de ``dimension de Kodaira" est 
surprenant, puisque Kodaira n'a jamais utilis\'e cette notion (sauf une fois, en 1975, alors 
qu'elle \'etait d\'ej\`a devenue d'usage courant), et ceci m\^eme dans sa 
classification des surfaces, qu'il base 
syst\'ematiquement sur l'\'etude des invariants num\'eriques de leurs mod\`eles minimaux).

\

La dimension canonique est l'invariant fondamental de la classification. Le th\'eor\`eme de prolongement d'Hartogs 
montre que c'est un invariant bim\'eromorphe de $X$. (Attention! 
Ce ne serait pas le cas si l'on avait consid\'er\'e 
$-K_X$). On peut donc, par invariance bim\'eromorphe, d\'efinir $\kappa(X)$ pour un espace analytique compact 
irr\'eductible $X$ quelconque (en en consid\'erant un mod\`ele lisse). 

Si $f:X\rightarrow Y$ est m\'eromorphe surjective, avec $X,Y$ de m\^eme dimension, 
on a: $\kappa(X)\geq \kappa(Y)$. En particulier: $X$ est de type g\'en\'eral si $Y$  l'est.

\

{\bf 1.7. Exemples:} 

\

{\bf 1.  Les courbes.} On a $3$ cas:

$\kappa(X)=-\infty$ ssi $X\cong \Bbb P^1$ ssi $g(X)=0$ ssi $K_X<0$.

$\kappa(X)=0$ ssi $X\cong \Bbb C/\Lambda $ ssi $g(X)=1$ ssi $K_X=0$.

$\kappa(X)=-\infty$ ssi $X\cong \Bbb D/\Gamma$ ssi $g(X)\geq 2$ ssi $K_X>0$.

On voit donc que pour les courbes, $\kappa$ \'equivaut \`a la ``g\'eom\'etrie" de $X$.

\

{\bf 2. Produits.} Si $X=T\times S$ est le produit de $2$ vari\'et\'es $T$ et $S$, alors 
on montre facilement que, de mani\`ere naturelle: $H^0(X,mK_X)=H^0(T,mK_T)\otimes H^0(S,mK_S),\forall m$, 
et donc que: 
$\kappa(X)=\kappa(T)+\kappa(S)$. 

Par exemple: Soit $T$ une vari\'et\'e complexe compacte et connexe. Alors: $\kappa(T\times \Bbb P^1)=-\infty$. 

Si $\kappa(S)=0$ (par exemple: si $S$ est une courbe elliptique), alors: $\kappa(T\times S)=\kappa(T)$. 
(On obtient ainsi des exemples simples de vari\'et\'es $X$ 
de dimension $n$, et de dimension canonique $\kappa\in \lbrace -\infty,0,1,...,n\rbrace$ arbitraire).

\

{\bf 3. Fibrations.} Soit $f:X\rightarrow Y$, 
une fibration. Alors on a: $ \kappa(Y)+\kappa(X_y)\leq\kappa(X)\leq dim(Y)+\kappa(X_y)$. 

L'in\'egalit\'e de droite est l'additivit\'e facile de {\bf 1.4.4}, mais celle de gauche est 
la conjecture $C_{n,m}$ d'Iitaka (l'hypoth\`ese $X$ K\" ahler est essentielle). 
 Cette conjecture est d\'emontr\'ee (Kawamata, Viehweg) lorsque 
$Y$ est de type g\'en\'eral (ie: $\kappa(Y)=dim(Y)>0$), auquel cas on a donc 
\'egalit\'e: $\kappa(X)= dim(Y)+\kappa(X_y)$.

\

Si $\kappa(X)\geq 0$, on d\'efinit par application de {\bf 1.4.4} \`a $K_X$ une (unique) fibration $J_X:X\to J(X)$, 
appel\'ee {\bf fibration d'Iitaka-Moishezon}, telle que: 
$dim(J(X))=\kappa(X)$, et $\kappa(X_j)=0$, si 
$X_j$ est sa fibre g\'en\'erale. Cette application d\'etermine donc la
 ``composante $\kappa=0$" de $X$ si $\kappa(X)\geq 0$. Observer 
cependant qu'il n'y a pas de relation simple entre $K_X$ et $J_X^*(K_{J(X)})$ en g\'en\'eral. 
En particulier, $J(X)$ n'est pas toujours de type g\'en\'eral, et $J_X$ ne fournit alors pas de 
``composante" de type g\'en\'eral de $X$.

(Remarquons que la fibration d'Iitaka-Moishezon, introduite dans [Ii] et [Mo], est appel\'ee simplement 
``fibration d'Iitaka", d'habitude. A nouveau, [Mo] semble \^etre pass\'e inaper\c cu).

\

{\bf 4. Surfaces.} On a ici $4$ classes de surfaces:

$\kappa(X)=-\infty$ ssi $X$ est bim\'eromorphe \`a un produit $\Bbb P^1\times C$, o\`u $C$ est une courbe.

$\kappa(X)=0$ ssi un rev\^etement \'etale fini de $X$ est bim\'eromorphe \`a un tore complexe (de dimension $2$), 
ou \`a une surface $K3$ (d\'eformation d'une quartique lisse de $\Bbb P^3$).

$\kappa(X)=1$ ssi $X$ est bim\'eromorphe \`a une fibration ``elliptique" $J_X=f:X\rightarrow C=J(X)$ 
sur une courbe $C$ telle que $mK_X=f^*(L)$, avec $L$ ample sur $C$, pour un $m>0$. (Dire que $f$ est elliptique 
signife que ses fibres lisses sont des courbes elliptiques).

$\kappa(X)=2$, alors $X$ est de ``type g\'en\'eral" (il n'y a pas de principe de classification connu; 
c'est dans cette classe que se trouve 
l'immense majorit\'e des vari\'et\'es de cette dimension (d'o\`u le terme, d\^u \`a Moishezon)).

\
{\bf 5. Vari\'et\'es de type g\'en\'eral.} Si $f:X\to Y$ est une fibration holomorphe et si $X$ est de type 
g\'en\'eral, alors $X_y$ est aussi de type g\'en\'eral, par additivit\'e facile {\bf 1.4.5} et $K_{X_y}=K_{X\vert X_y}$. Par restriction \`a un sous-espace ad\'equat de $Z$, on en d\'eduit que le membre g\'en\'eral d'une famille de sous-vari\'et\'es $(V_z)_{z\in Z}$ est aussi de type g\'en\'eral si ces sous-vari\'et\'es recouvrent $X$.

\

\

{\bf C. Vari\'et\'es unir\'egl\'ees.} 

\

{\bf 1.8. D\'efinition:} On dit que $X$ (K\" ahler compacte connexe) est {\bf unir\'egl\'ee} s'il 
existe une application m\'eromorphe surjective $f:T\times \Bbb P^1\rightarrow X$ 
(dans laquelle $T$ d\'epend de $X$, et la restriction de $f$ \`a $t\times \Bbb P^1$ est non-constante, pour $t\in T$ g\'en\'erique).

\

On montre, \`a l'aide de la th\'eorie de la vari\'et\'e de Chow (qui d\'ecrit en particulier 
les d\'eformations de $f$ par des 
espaces de param\`etres analytiques), que cette condition signifie aussi que $X$ 
est recouverte par des {\bf courbes rationnelles} (ie: des images d'applications holomorphes non-constantes
$r:\Bbb P^1\rightarrow X$), et que l'on peut choisir $f:T\times \Bbb P^1\rightarrow X$ g\'en\'eriquement finie.

\

On en d\'eduit la:

\

{\bf 1.9. Proposition:} Si $X$ est unir\'egl\'ee, alors: $\kappa(X)=-\infty$.

\

L'une des conjectures fondamentales de la classification est la r\'eciproque, \'etablie en dimension 
$3$ si $X$ est projective ([Mi], [M]), mais seulement en dimension $2$ pour $X$ K\" ahler:

\

{\bf 1.10. Conjecture ``$-\infty$":} Si $\kappa(X)=-\infty$, alors $X$ est unir\'egl\'ee.

\

Cette conjecture n\'ecessite l'hypoth\`ese $X$ K\" ahler (comme le montrent les surfaces de Hopf).

\

{\bf 2. LES TROIS G\'EOM\'ETRIES ``PURES".}

\

{\bf A. Les trois g\'eom\'etries pures.}

\

Elles sont d\'efinies en dimension $n$ 
par les trois classes suivantes de vari\'et\'es $X$ (qui g\'en\'eralisent respectivement 
les courbes hyperboliques, plates et sph\'eriques):

\

{\bf 2.1. Vari\'et\'es de type g\'en\'eral:} Ce sont les $X$ telles que $\kappa(X)=n$. Les exemples 
les plus simples sont les quotients de domaines born\'es sym\'etriques, et les hypersurfaces lisses 
de degr\'e au moins $n+3$ de $\Bbb P^{n+1}$.

\

{\bf 2.2. Vari\'et\'es avec $\kappa(X)=0$.} Les exemples les plus simples sont les tores complexes 
($\Bbb C^n/\Lambda)$, $\Lambda$ r\'eseau de $\Bbb C^n$, les vari\'et\'es hyperk\"ahl\'eriennes, 
et les vari\'et\'es de Calabi-Yau. Les hypersurfaces lisses de degr\'e $n+2$ dans $\Bbb P^{n+1}$ 
appartiennent \`a cette classe. Pour $n=1$, on obtient les courbes elliptiques, pour $n=2$, des vari\'et\'es 
hyperk\"ahl\'eriennes: les surfaces $K3$, et pour $n\geq 3$, des vari\'et\'es de Calabi-Yau, dont 
les quintiques de $\Bbb P^4$ fournissent les premiers exemples).

\

{\bf 2.3. Les vari\'et\'es $\kappa$-Rationnellement Connexes.} (Ou ``$\kappa$-RC'', en abr\'eg\'e). Ce sont les $X$ telles que
, pour toute fibration (m\'eromorphe) $f:X\to Y$, on ait: $\kappa(Y)=-\infty$. (ie: $X$ ne  ``domine" que des 
vari\'et\'es $Y$ telles que $\kappa(Y)=-\infty$).

Si l'on prend $f:=id_X:X\to X$, on a donc: $\kappa(X)=-\infty$ si $X$ est $\kappa$-RC. 
La conjecture $-\infty$ affirme donc que $X$ est unir\'egl\'ee si $X$ est $\kappa$-RC. 

Mais
la condition ``unir\'egl\'ee" est beaucoup plus faible que $\kappa$-RC, car si $f:X\to Y$ est une 
fibration, alors $Y$ est $\kappa$-RC si $X$ l'est. Tandis que $X=\Bbb P^1\times T$ 
est unir\'egl\'ee, pour $T$ arbitraire, mais $\kappa$-RC ssi $T$ l'est.
Si $C$ est une courbe elliptique, par exemple, $X=\Bbb P^1\times C$ est unir\'egl\'ee, mais pas $\kappa$-RC.

\

La propri\'et\'e $\kappa$-RC peut n\'emmoins \^etre d\'ecrite g\'eom\'etriquement, 
si l'on admet la conjecture $-\infty$.

\

{\bf 2.4. D\'efinition:} La vari\'et\'e $X$ est {\bf rationnellement connexe} (RC en abr\'eg\'e) si deux 
points g\'en\'eriques de $X$ peuvent \^etre joints par une courbe rationnelle de $X$.

\

{\bf 2.5. Remarques.} Si $X$ est lisse, il suffit pour \^etre RC, que $2$ points g\'en\'eriques de $X$ soient 
contenus dans une r\'eunion finie {\bf connexe} de courbes rationnelles de $X$. La lissit\'e est essentielle 
(prendre un c\^one sur une courbe elliptique pour le voir). De plus, si $X$ est RC, tout ensemble fini de $X$ 
(et non seulement les ensembles \`a $2$ \'el\'ements) est contenu dans une courbe rationnelle de $X$. 
Ces r\'esultats difficiles sont d\'emontr\'es dans [Ko-Mi-Mo]. Voir [De97] pour un expos\'e tr\`es clair.

\

{\bf 2.6. Exemples:} $\Bbb P^n$, les Grassmanniennes, et plus g\'en\'eralement, 
les vari\'et\'es rationnelles ou unirationnelles sont RC. Les vari\'et\'es de Fano ($-K_X$ ample) le sont 
aussi, et donc aussi les hypersurfaces lisses de degr\'e au plus $n+1$ de $\Bbb P^{n+1}$, 
par la formule d'adjonction.

\

Une vari\'et\'e $X$ qui est RC est unir\'egl\'ee, et a donc $\kappa(X)=-\infty$. Par ailleurs, 
si $f:X\to Y$ est surjective, il est imm\'ediat que $Y$ est \'egalement RC, donc $\kappa(Y)=-\infty$. Donc 
$X$ est $\kappa$-RC si $X$ est RC. R\'eciproquement, on montre facilement (\`a l'aide du quotient rationnel 
d\'ecrit ci-dessous et de [G-H-S] en particulier) la:

\

{\bf 2.7. Proposition:} Si la conjecture $-\infty$ est vraie, alors $X$ est RC ssi $X$ est $\kappa$-RC.

(Voir {\bf 2.12} ci-dessous pour les ingr\'edients de la d\'emonstration).

\

{\bf B. Le quotient rationnel.} 

\

On a d\'efini en {\bf 1.7.3} la fibration d'Iitaka-Moishezon de $X$ si $\kappa(X)=0$. 
Cette fibration d\'etermine dans un sens satisfaisant la ``composante $\kappa(X)=0$"  dans ce cas.

On peut aussi construire la ``composante" RC de $X$ (qui est donc aussi sa composante pure $\kappa$-RC si 
la conjecture $-\infty$ est vraie).

\

{\bf 2.8. Th\'eor\`eme:} Il existe une (unique) fibration $r_X:X\to R(X)$ telle que:

{\bf 1.} Les fibres de $r_X$ sont RC.

{\bf 2.} $R(X)$ n'est pas unir\'egl\'e.

On appelle $r_X$ le {\bf quotient rationnel} de $X$.

\ 

{\bf 2.9. Remarques:} 

{\bf 1.} La conjecture $-\infty$ affirme que $\kappa(R(X))\geq 0$, et donc que 
l'on peut d\'efinir la compos\'ee $J_{R(X)}\circ r_X:X\to J(R(X))$ dans tous les cas.

{\bf 2.} On peut d\'efinir de mani\`ere analogue un {\bf quotient $\kappa$-RC}, en 
rempla\c cant RC par $\kappa$-RC dans la condition {\bf 3.7.1}. Ces deux 
quotients coincident, si {\bf 1.10} est vraie.

{\bf 3.} Observer que $dim(R(X))<dim(X)$ ssi $X$ est unir\'egl\'ee.

\

{\bf 2.10. indication} sur la d\'emonstration de {\bf 2.8}: On construit ([C81], [Ko-Mi-Mo]) 
$r_X$ comme le quotient de $X$ pour 
la relation d'\'equivalence sur $X$, lisse, qui identifie deux points lorsqu'ils sont 
contenus dans une r\'eunion 
connexe de courbes rationnelles. Les fibres g\'en\'eriques de $r_X$ sont alors lisses, donc RC par {\bf 2.5}.
On conclut enfin que $R(X)$ n'est pas unir\'egl\'e \`a l'aide du:

\

{\bf 2.11. Th\'eor\`eme ([G-H-S]):} Soit $f:X\to Y$ une fibration telle que $X_y$ soit RC. Alors:

1. Si $Y$ est une courbe, $f$ admet une section. En particulier, $X$ est RC si $Y=\Bbb P^1$.

2. Si $Y$ est $RC$, $X$ aussi.

\

(La d\'emonstration de {\bf 2.11} (voir [De02] pour un expos\'e accessible) est donn\'ee dans 
le cas o\`u $X$ est projective, mais s'applique au cas K\" ahler par 
certains r\'esultats de [C81]).

\

{\bf 2.12. indication} sur la d\'emonstration de {\bf 2.7}: Si $R(X)$ n'est pas 
un point, alors $\kappa(R(X))\geq 0$, par la conjecture $-\infty$. Contradiction si $X$ est $\kappa$-RC. 
Donc $R(X)$ est un point, et $X$ est RC.

\

{\bf C. Pourquoi les structures orbifoldes.}

\

{\bf 2.13. Une id\'ee naturelle} pour construire les ``composantes" $\kappa=0$ et $\kappa$-RC de $X$ 
consiste \`a it\'erer l'application  $(J_{R(X)}\circ r_X)$ pr\'ec\'edente; 
c'est-\`a-dire \`a l'appliquer \`a $X$, puis \`a $J(R(X))=(J\circ r)$, 
$J(R(J(R(X))))=(J\circ r)^2$, etc... Cette suite de fibrations stationne quand la base $B$, obtenue apr\`es 
$k$ it\'erations $c:=(J\circ r)^k, k\leq n$, satisfait: $B=R(B)=J(R(B))$, ce qui ne se produit clairement 
que si $B$ est 
une vari\'et\'e de type g\'en\'eral, ou un point. On obtiendrait ainsi, par \'elimination des 
``composantes" successives des types $\kappa=0$ et $\kappa$-RC de $X$, une base 
qui serait la ``composante" de type g\'en\'eral de $X$ (\'eventuellement triviale, 
c'est-\`a-dire r\'eduite \`a un point).

Remarquons que l'on peut envisager de construire directement ``par en bas", la fibration 
$c:=(J\circ r)^n$ compos\'ee de toutes les fibrations pr\'ec\'edentes en 
construisant une fibration $f:X\to Y$ de base $Y$ de type g\'en\'eral, et 
maximale pour cette propri\'et\'e. On peut facilement construire une telle $f$ , mais la construction 
n'est pas stable par rev\^etement \'etale fini. Et ne peut en fait pas \^etre corrig\'ee en 
prenant en consid\'eration ces rev\^etements. ([C01],[B-T]). 
 
Nous allons voir sur un exemple simple ({\bf 2.14} ci-dessous) que des ``composantes" de type g\'en\'eral 
peuvent exister dans $X$, qui ne sont pas ``r\'ev\'el\'ees" par une telle it\'eration, 
\`a cause de la pr\'esence de fibres multiples. Cet exemple montre aussi que la dimension des fibres de la 
fibration $c$ n'est pas 
invariante par rev\^etement \'etale fini.

\ 

{\bf 2.14. Un exemple justificatif.} Soient $E$ (resp. $H$) une courbe elliptique (resp. hyperelliptique de genre $g\geq 2$)
, et $t$ (resp. $\vartheta$) une translation 
d'ordre $2$ sur $E$ (resp. l'automorphisme hyperelliptique de $H$). Soit $X':=E\times H$, et 
$a:=(t\times \vartheta: X'\to X'$ l'involution (sans point fixe) obtenue par action diagonale. 

On note $X:=X'/<a>$ le quotient de $X'$ par $a$, et par $u:X'\to X$ et $v:H\to H/<\vartheta>\cong \Bbb P^1$ 
les applications quotient naturelles.

Alors $J_X:X\to H/<\vartheta>=J(X)$ est la projection naturelle. Elle a pour base $\Bbb P^1$, 
et fibre g\'en\'erique $E$. N\'eammoins, $X$ a certainement une ``composante" de type g\'en\'eral 
non triviale (r\'ev\'el\'ee par la pseudom\'etrique de Kobayashi, ou le groupe fondamental, si l'on veut 
qu'une telle notion soit invariante par rev\^etement \'etale). Pour cettte raison, la base $\Bbb P^1$ de 
la fibration $J_X$ doit \^etre consid\'er\'ee comme \'etant de type g\'en\'eral, si l'on veut 
avoir une construction 
compatible avec les rev\^etements \'etales.

\

$$\xymatrix{X'=E\times H\ar[r]^u\ar[d]_{J_{X'}} & X\ar[d]^{J_X} \\
H\ar[r]^v & \Bbb P^1
}$$

\
Une donn\'ee g\'eom\'etrique a n\'eammoins gard\'e la trace de la construction pr\'ec\'edente: les fibres multiples 
de $J_X$. Celles-ci sont de multiplicit\'e $2$, situ\'ees au-dessus des images $b_j, j= 1,...,2g+2$, par 
$v$ des points hyperelliptiques de $H$. On va donc d\'efinir le fibr\'e canonique de la base $B$ de $J_X$ par:
 $K_B:=K_{\Bbb P^1}+\sum_{j=1,...,2g+2}(1-1/2).b_j$, de telle sorte que 
son image r\'eciproque par $v$ soit $K_H$.

\

{\bf 2.15. La version orbifolde} Nous allons maintenant g\'en\'eraliser cette construction, 
en introduisant {\bf (\S 4)} 
une structure orbifolde 
sur la base d'une fibration $f:X\to Y$ arbitraire, et montrer que les structures 
orbifoldes ainsi obtenues permettent 
de r\'esoudre naturellement le probl\`eme de d\'ecomposition pos\'e initialement, en consid\'erant les 
versions orbifoldes des fibrations $r,J,c$ \'evoqu\'ees en {\bf 2.9} ci-dessus. 

L'id\'ee sous-jacente est que, tout comme dans l'exemple {\bf (2.14)} ci-dessus, il 
existe un rev\^etement ramifi\'e fini $Y'$ de $Y$ tel que par ce changement de base les fibres de $f$ perdent 
leurs multiplicit\'es, et que l'espace $X'$ d\'eduit de $X$ soit non-ramifi\'e sur $X$. 
Un tel changement de base n'existe pas, en g\'en\'eral, et la structure orbifolde que l'on introduit sur $Y$ 
est un substitut virtuel pour un tel $Y'$. On va justement v\'erifier que les orbifoldes ainsi introduites 
se comportent exactement comme dans le cas ``classique" ou un tel rev\^etement existerait.

Cette notion de base orbifolde nous permettra de d\'efinir {\bf (\S 6.A)} 
les vari\'et\'es sp\'eciales comme \'etant celles 
n'ayant pas de fibration non-constante ayant une base orbifolde de type g\'en\'eral. Cette classe 
contient les vari\'et\'es de type RC ou $\kappa=0$.

On construira directement (par ``en bas"), {\bf (\S 7.A}) sur toute $X$ une unique fibration (le ``coeur") 
\`a fibres sp\'eciales et base 
orbifolde de type g\'en\'eral (ou un point ssi $X$ est sp\'eciale). 

L'it\'eration de la suite de fibrations $(J\circ r)$, prise au sens orbifolde, 
nous conduira en {\bf 7.B} \`a la d\'ecomposition canonique du ``coeur".

\

{\bf 3. ORBIFOLDES.}

\

{\bf 3.1.} Soit $Y$ une vari\'et\'e complexe compacte et connexe. Une structure 
{\bf d'orbifolde} sur $Y$ est la donn\'ee 
d'un $\Bbb Q$-diviseur $\Delta:=\sum_{j\in J}a_j.D_j$, avec $0<a_j\leq 1$, 
$a_j\in \Bbb Q,\forall j$, et $D_j$ 
des diviseurs irr\'eductibles distincts de $Y$. (On dira que $\Delta$ est \`a 
multiplicit\'es ``standard" si $a_j=(1-1/m_j), 
m_j\in \Bbb Z, \forall j\in J$).

On notera $(Y/\Delta)$ l'orbifolde ainsi d\'efinie.

\

Son {\bf fibr\'e canonique} est le $\Bbb Q$-diviseur $K_Y+\Delta$ sur $Y$.

\

Sa {\bf dimension de Kodaira} est: 
$\kappa(Y/\Delta):=\kappa(Y,K_Y+\Delta):=\kappa(Y,m_0.(K_Y+\Delta)$, o\`u l'on a choisi
 $m_0\in \Bbb Z^+$ tel que $m_0.a_j\in \Bbb Z, \forall j$.

On a, bien s\^ur, toujours: $\kappa(Y/\Delta)\geq \kappa(Y)$, et $\kappa(Y/\Delta)\in \lbrace -\infty,0,1,...,n\rbrace$. 

\

{\bf 3.2. Exemple:} Si $Y=\Bbb P^p$, si $d_j$ est le degr\'e de $D_j$, et si 
$\delta:=deg(K_{(\Bbb P^p/\Delta)}):=-(p+1)+\sum_{j\in J}a_j.d_j$, alors: 
$\kappa(\Bbb P^p/\Delta)=-\infty$ (resp. $0$; resp. $p$) ssi: 
$\delta <0$ (resp. $\delta=0$; resp. $\delta>0$). 

\

{\bf 3.3. Fibration d'Iitaka-Moishezon.} Elle peut \^etre d\'efinie par application de 
{\bf 1.4.4.} \`a son fibr\'e canonique, avec les m\^emes propri\'et\'es que si $\Delta=\emptyset$,
 pour toute orbifolde $(Y/\Delta)$ si $\kappa(Y/\Delta)\geq 0$.

\

{\bf 4. LA BASE ORBIFOLDE D'UNE FIBRATION.}

\

{\bf A. G\'en\'eralit\'es.}

\

{\bf 4.0. Multiplicit\'es.} Soit $f:X\rightarrow Y$ une fibration holomorphe, avec $X,Y$ lisses. 

Un diviseur effectif $E$ de $X$ est dit {\bf $f$-exceptionnel} si $codim_Y(f(E))\geq 2$.

Soit $D\subset Y$ un diviseur irr\'eductible, et $f^*(D)=\sum_{j\in J}(m_jD_j)+E$, 
o\`u $E$ est $f$-exceptionnel, et o\`u les $D_j$ sont les composantes irr\'eductibles 
non $f$-exceptionnelles de $f^{-1}(D)$. Donc: $f(D_j)=D, \forall j$.

On d\'efinit alors: $m(f,D):=inf_{j\in J}\lbrace m_j\rbrace$. C'est la {\bf multiplicit\'e de la 
fibre de $f$ 
au-dessus du point g\'en\'erique de $D$}.

Remarquons que l'on a donc: $m(f,D)=1$, sauf pour un nombre fini de diviseurs $D$, contenus dans le 
lieu de $Y$ au-dessus duquel $f$ n'est pas submersive.

\

{\bf 4.1. Remarque:} La multiplicit\'e classique est d\'efinie diff\'eremment,
 par: $m^+(f,D):=pgcd_{j\in J}\lbrace m_j\rbrace$.

\

{\bf 4.2. D\'efinition:} La {\bf base orbifolde de $f$} est: $(Y/\Delta(f))$, o\`u :
$\Delta(f):=\sum_{D\subset Y}(1-1/m(f,D)).D$. (Cette somme est finie,
 puisque le coefficient de $D$ s'annule ssi $m(f,D)=1$).

\

{\bf 4.3. Exemples:}

1. Si $X_y$ est RC et si $Y$ est projective, on d\'eduit de {\bf 2.11} que $\Delta(f)=\emptyset$.

2. Si $X_y$ est une courbe elliptique, on peut avoir $\Delta(f)\neq \emptyset$ (voir l'exemple {\bf 2.14}).

On d\'efinit enfin: $\kappa(Y,f):=\kappa(Y/\Delta(f)):=\kappa(Y,K_Y+\Delta(f))$.

\

{\bf 4.4. Remarque:} Cette dimension n'est pas un invariant bim\'eromorphe de $f$. 
C'est-\`a-dire que l'on n'a 
pas n\'ecessairement $\kappa(Y,f)=\kappa(Y',f')$ si $f$ et $f':X'\to Y'$ 
sont \'equivalentes (not\'e: $f\sim f'$)
, c'est-\`a-dire s'il existe un 
diagramme commutatif $u:X'\to X$, $v:Y'\to Y$ avec $f\circ u=v\circ f':X'\to Y$, 
dans lequel $u$ et $v$ sont bim\'eromorphes.

\

$$\xymatrix{X'\ar[r]^u\ar[d]_{f'} & X\ar[d]^f \\
Y'\ar[r]^v & Y
}$$

\

Dans une telle situation, on a alors:

{\bf 1.} $v_*(\Delta(f')=\Delta(f)$, et donc:

{\bf 2.} $\kappa(Y',f')\leq \kappa(Y,f)$.

Mais on peut avoir in\'egalit\'e stricte (seulement si $\kappa(Y)=-\infty$). D'o\`u:

\

{\bf 4.5. D\'efinition:} $\kappa(f):=inf \lbrace \kappa(f'); {f'\sim f}\rbrace$.

On dit que $f$ est {\bf de type g\'en\'eral} si $\kappa(f)=dim(Y)>0$.

\

(Si $Y$ est de type g\'en\'eral, $f$ est donc de type g\'en\'eral). Mais il existe des 
fibrations de type g\'en\'eral de base $\Bbb P^p$, par exemple celle d\'efinie en {\bf 2.14}).

\

On a un crit\`ere pour tester si $\kappa(Y,f)=\kappa(f)$.

\

{\bf 4.6. D\'efinition:} $f':X'\rightarrow Y'$ est {\bf nette} s'il existe $u:X'\rightarrow X$, 
avec $X$ lisse, telle que 
tout diviseur $E'\subset X'$ qui est $f'$-exceptionnel est $u$-exceptionnel.

\

Une fibration $f$ a toujours des repr\'esentants nets, obtenus par aplatissement de $f$ par 
un changement de base bim\'eromorphe $v:Y'\rightarrow Y$, et d\'esingularisation $d:X'\to X$ de
 l'espace $X"$ obtenu par ce changement de base.

\ 

$$\xymatrix{X'\ar[r]^d \ar[rd]_{f'} & X"\ar[r]^{u"}\ar[d]_{f"} & X\ar[d]^f \\
& Y'\ar[r]^v & Y
}$$

\

{\bf 4.7. Proposition:} $\kappa(Y,f)=\kappa(f)$ si $f$ est nette.

\

{\bf 4.8. Remarque:} On d\'emontre la proposition ci-dessus par une seconde d\'efinition directe de 
$\kappa(f):=\kappa(X,F_f)$, o\`u $F_f\subset \Omega _X^p$ est le faisceau coh\'erent de rang $1$ 
obtenu par saturation de $f^*(K_Y)$ dans $\Omega _X^p$.

Les fibrations de type g\'en\'eral sur $X$ lisse sont alors naturellement en bijection avec les sous-faisceaux 
{\bf de Bogomolov} de $X$, qui sont les sous-faisceaux satur\'es coh\'erents $F$ de rang $1$ de $\Omega_X^p$, 
pour un $p>0$, tels que $\kappa(F)=p$, c'est-\`a-dire maximale, par l'in\'egalit\'e classique de Bogomolov. 

\

{\bf 4.9. Remarque:} Si $f:X\rightarrow Y$ est une fibration, et $E:=\sum_{h\in H} (1-1/n_h).E_h$ 
un diviseur orbifolde ``standard" sur $X$, 
alors la restriction $E_y$ de $E$ \`a $X_y$ est un diviseur orbifolde standard sur $X_y$. On d\'efinit
 de mani\`ere analogue au cas ci-dessus, o\`u $E=\emptyset$,
la base orbifolde $\Delta(f,E):=\sum_{D\subset Y} (1-1/m(f,E;D)).D$ de $(f,E)$, et la dimension de Kodaira $\kappa(f,E)$,  de telle sorte que si 
$g:Z\rightarrow X$ et $f$ sont nettes, alors: $\Delta(f\circ g)=\Delta(f,\Delta(g))$.

(Plus pr\'ecis\'ement: $m(f,E;D)=inf_{j\in J}\lbrace m_j.n_j\rbrace$, o\`u:
  $n_j=n_h$ si $Dj=E_h$, $h\in H$, et $n_j=1$ sinon).

\

$$\xymatrix{Z\ar[r]^g\ar[rd]_{f\circ g} & (X/\Delta (g))\ar[d]^f \\
& Y
}$$

\

{\bf 5. ADDITIVIT\'E ORBIFOLDE.}

\

On \'etend au cadre orbifolde la conjecture $C_{n,m}$ d'Iitaka.

\

{\bf 5.1. Conjecture $C_{n,m}^{orb}$:} Soit $f:(X/E)\rightarrow Y$ une fibration 
\`a croisements normaux (dans un sens
ad\'equat naturel), $E$ 
\'etant un diviseur orbifolde sur $X$. Alors: $\kappa(X/E)\geq \kappa (X_y/E_y)+\kappa(Y/(\Delta(f,E))$.

\

L'un des outils essentiels des pr\'esentes consid\'erations est le:

\

{\bf 5.2. Th\'eor\`eme:} Soit $f:(X/E)\rightarrow Y$ comme ci-dessus. 

Si $(f,E)$ est de type g\'en\'eral (ie: 
si $\kappa(Y/\Delta(f,E))=dim(Y)$), alors: $\kappa(X/E)=dim(Y)+\kappa(X_y/E_y)$.

\

La d\'emonstration est une adaptation au cas orbifolde de celle, classique, de Viehweg. Bien que 
les techniques de d\'emonstrations soient les m\^emes, le cadre orbifolde en \'etend consid\'erablement le 
domaine d'application quand $\kappa(X_y)=-\infty$.

\

Dans le cas particulier o\`u $E=\emptyset$, on obtient:

\

{\bf 5.3. Corollaire 1.} Soit $f:X\rightarrow Y$ une fibration de type g\'en\'eral 
(ie: $\kappa(f)=dim(Y)>0$). Alors: $\kappa(X)=\kappa(X_y)+dim(Y)$.

\

Si $\kappa(X)=0$, on a donc:

\

{\bf 5.4. Corollaire 2.} Si $\kappa(X)=0$, il n'existe pas de fibration de type g\'en\'eral $f:X\rightarrow Y$. 
(On dira que $X$ est ``sp\'eciale").

\

{\bf 5.5. Corollaire 3.} Soit $f:X\rightarrow Y$ et $g:Z\rightarrow X$ des fibrations. On suppose que:

{\bf 1.} $(f\circ g)$ est de type g\'en\'eral.

{\bf 2.} La restriction $g_y:Z_y\rightarrow X_y$ est de type g\'en\'eral (pour $y\in Y$ g\'en\'eral).

Alors: $g:Z\rightarrow X$ est de type g\'en\'eral.

\

Ce r\'esultat joue un r\^ole crucial dans la suite.

\

{\bf id\'ee de la d\'emonstration de 5.3:} On suppose que $dim(Y)=1$ et 
que $m^+(f,D)=m(f,D), \forall D\subset Y$. On a alors un diagramme commutatif dans lequel 
$u$ est \'etale fini, et $v$ fini, avec: $v^*(K_Y+\Delta(f))=K_{Y'}$.

\

$$\xymatrix{X'\ar[r]^{u=etale}\ar[d]_{f'} & X\ar[d]^f \\
Y'\ar[r]^{v} & Y
}$$

\
On se ram\`ene au cas classique en observant que:

. $v^*(m(K_Y+\Delta(f)))=m.K_{Y'}$, et donc (par platitude de $f$), que:

. $v^*(f_*(m.K_{(X/(Y/\Delta(f)))}))=(f')_*(m.K_{X'/Y'})$, la notation $K_{X/Y}$ d\'esignant 
comme d'habitude le fibr\'e canonique relatif (y compris si $Y$ est une orbifolde).

Le cas g\'en\'eral se traite de la m\^eme fa\c con, avec des d\'etails techniques additionnels. 

\

{\bf 6.  VARI\'ET\'ES SP\'ECIALES.}

\

{\bf A. G\'en\'eralit\'es.}

\

{\bf 6.1. D\'efinition:} On dit que $X$ est {\bf sp\'eciale} s'il n'existe pas de fibration (m\'eromorphe) 
$f:X\rightarrow Y$ de type g\'en\'eral.

\

{\bf 6.2. Exemples:} 

\

{\bf 0.} Vari\'et\'es de type g\'en\'eral: Si $X$ est de type g\'en\'eral, alors $X$ n'est pas 
sp\'eciale (car $id_X:X\rightarrow X$ est alors une fibration de type g\'en\'eral).

\

{\bf 1.} Courbes: La courbe $X$ est sp\'eciale ssi elliptique ou $\Bbb P^1$, car 
une fibration d\'efinie sur $X$ est soit constante, soit $id_X$, l'identit\'e de $X$.

\

{\bf 2.} Rev\^etements: Si $g:X\rightarrow X'$ est une application m\'eromorphe surjective, et si $X$ est 
sp\'eciale, alors $X'$ aussi. 

Si $g$ est un rev\^etement \'etale fini, et si $X'$ est sp\'eciale, alors 
$X$ est aussi sp\'eciale. 

La preuve de cette derni\`ere assertion d'apparence triviale est difficile, et utilise le Corollaire 5.5 
pr\'ec\'edent. 

\

{\bf 3.} Surfaces: Une surface $X$ est sp\'eciale ssi elle est dans l'une des $3$ classes suivantes:

. $\kappa(X)=-\infty$: $X$ est bim\'eromorphe \`a $\Bbb P^1\times C$, 
avec $C$ une courbe elliptique ou $\Bbb P^1$.

. $\kappa(X)=0$: $X$ a un rev\^etement \'etale fini bim\'eromorphe \`a $\Bbb C^2/\Lambda$ (tore complexe), ou 
\`a une surface $K3$.

. $\kappa(X)=1$: $X$ a un rev\^etement \'etale fini $X'$ qui est une fibration elliptique 
$f':X'\rightarrow C$, de base une courbe $C$ elliptique ou $\Bbb P^1$, $f'$ ayant au plus $2$
 fibres multiples si $C\cong \Bbb P^1$, et aucune si $C$ est elliptique.

\

Plus g\'en\'eralement, il existe pour tout $n>0$ et tout $\kappa\in \lbrace -\infty, 0,1,...,(n-1)\rbrace$, 
des vari\'et\'es sp\'eciales $X$ de dimension $n$ et telles que $\kappa(X)=\kappa$.

\

{\bf 4.} Fibrations \`a fibres et base sp\'eciales: Soit $f:X\to Y$ une fibration 
telle que $X_y$ (la fibre g\'en\'erale) et $Y$ sont sp\'eciales. 
Alors: $X$ n'est pas sp\'eciale, en g\'en\'eral (voir l'exemple justificatif 2.14). Cependant:

\

{\bf 6.3. Th\'eor\`eme:} Si $f:X\to Y$ est une fibration 
telle que $X_y$ (la fibre g\'en\'erale) et $Y$ sont sp\'eciales, et si, de plus, 
$f$ n'a pas de fibre multiple en codimension $1$ (par exemple: si $f$ a une section)
, alors $X$ est sp\'eciale.

\

La d\'emonstration de {\bf 6.3} repose sur le: 

\

{\bf 6.4. Th\'eor\`eme:} Soit $f:X\rightarrow Y$ et $g:X\to Z$ des fibrations telles que: $g$ 
est de type g\'en\'eral, et $X_y$ (la fibre g\'en\'erale de $f$) est sp\'eciale. Alors: il existe 
une (unique) fibration $h:Y\to Z$ telle que $g=h\circ g$. (Les fibrations \`a fibres 
sp\'eciales ``dominent" les fibrations de type g\'en\'eral). 

\

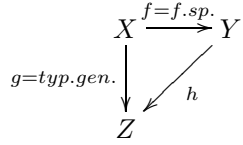
$$\xymatrix{X\ar[r]^{f=f.sp.}\ar[d]_{g=typ.gen.} & Y\ar[ld]^h \\
Z 
}$$

\

{\bf id\'ee de la d\'emonstration:} Si $Z$ (et non seulement $g$) est de type g\'en\'eral, la famille 
$(V_y:=g(X_y))_{y\in Y}$ de sous-vari\'et\'es de $Z$ param\'etr\'ee par $Y$ recouvre $Z$. Il s'agit de montrer que 
$dim(V_y)=0$. Sinon, on d\'eduit de ({\bf 1.7.5}) que $V_y$ est de type g\'en\'eral.
Ceci contredit l'hypoth\`ese que $X_y$ est sp\'eciale. Le cas g\'en\'eral s'obtient par un argument analogue, 
apr\`es avoir montr\'e que (la factorisation de Stein de) 
la restriction de $g$ \`a $V_y$ est une fibration de type g\'en\'eral, pour $y\in Y$ g\'en\'erique.

\

{\bf B. Vari\'et\'es RC, ou avec $\kappa=0$.}

\

Il s'agit des exemples fondamentaux de vari\'et\'es sp\'eciales, qui g\'en\'eralisent le cas des courbes:

\

{\bf 6.5. Th\'eor\`eme:} $X$ est sp\'eciale dans les $2$ cas suivants:

{\bf 1.} $X$ est rationnellement connexe.

{\bf 2.} $\kappa(X)=0$.

\

Le cas {\bf 2.} est une simple reformulation du corollaire 5.4 ci-dessus. Le cas 
{\bf 1.} peut \^etre d\'eduit du r\'esultat suivant, appliqu\'e aux courbes rationnelles de $X$ (qui sont 
sp\'eciales):

\

{\bf 6.6. Th\'eor\`eme:} Soit $X$ une vari\'et\'e dans laquelle $2$ points g\'en\'eraux peuvent \^etre 
joints par une chaine connexe de sous-vari\'et\'es sp\'eciales. Alors: $X$ est sp\'eciale.

\

{\bf id\'ee de la d\'emonstration:} Sinon, soit $f:X\to Y$ de type g\'en\'eral. Les chaines de 
sous-vari\'et\'es sp\'eciales de $X$ qui rencontrent $X_y$, fibre g\'en\'erale de $f$, sont contenues dans $X_y$, 
puisque la restriction de $f$ \`a une sous-vari\'et\'e sp\'eciale de $X$ en position g\'en\'erale est 
encore de type g\'en\'eral. Ceci contredit le fait que $dim(Y)>0$, et que l'on peut joindre un point g\'en\'eral de $X_y$ 
\`a un point g\'en\'eral de $X$ par une telle chaine.

\

{\bf 6.7. Remarques:} 

1. On ne sait pas si les vari\'et\'es $\kappa$-RC sont sp\'eciales. C'est vrai cependant pour leur 
variante orbifolde (les vari\'et\'es avec $\kappa_+=-\infty$, d\'efinies en {\bf 7.5} ci-dessous).

2. On verra en {\bf 7.7} ci-dessous, que, r\'eciproquement, les vari\'et\'es sp\'eciales sont des tours de 
fibrations en orbifoldes avec: soit $\kappa_+=-\infty$, soit $\kappa=0$.

\

{\bf C. Kobayashi-Ochiai orbifolde.}

\

Un nouvel exemple de vari\'et\'e sp\'eciales est fourni par le:

\

{\bf 6.8. Th\'eor\`eme:} Soit $\varphi:\Bbb C^n\to X$ une application m\'eromorphe dominante (ie: de rang $n=dim(X)$ 
en un point au moins de $\Bbb C^n$). Alors: $X$ est sp\'eciale.

\

{\bf 6.9. Remarque:} 
On peut donner des versions plus g\'en\'erales de 6.8, avec la m\^eme conclusion, en y 
rempla\c cant $\Bbb C^n$ 
par $W\times \Bbb C$, $W$ quasi-projective arbitraire, pourvu que $\varphi(W\times\lbrace 0\rbrace)=a\in X$. On 
obtient ainsi une version transcendante du fait que les vari\'et\'es RC sont sp\'eciales..

\

On d\'eduit {\bf 6.8} de la version orbifolde suivante du th\'eor\`eme d'extension de 
Kobayashi-Ochiai [KO75]:

\

{\bf 6.10. Th\'eor\`eme:} Soit $U$ un ouvert de Zariski d'une vari\'et\'e complexe connexe $V$, et 
$\varphi:U\to X$ une application m\'eromorphe. Soit $f:X\to Y$ une fibration de type g\'en\'eral. On suppose 
que $\psi:=f\circ \varphi:U\to Y$ est dominante. Alors $\psi$ s'\'etend m\'eromorphiquement \`a $V$.

\

Le cas trait\'e par Kobayashi-Ochiai est celui o\`u $X=Y$ est une vari\'et\'e de type g\'en\'eral 
(Dans {\bf 6.13}, $Y$ 
peut \^etre un espace projectif, comme dans l'exemple {\bf 2.14}).

\

{\bf id\'ee de la d\'emonstration:}

$$\xymatrix{U:=(V-D)\ar[r]^>>>>{\varphi}\ar[rd]_{\psi=f\circ \varphi} & X\ar[d]^f \\
& Y}$$

La d\'emonstration de [KO75] (avec $X=Y$) repose sur le fait que $\varphi^*(mK_Y)\subset mK_U$. On 
peut l'adapter \`a notre situation orbifolde en montrant que 
$f^*(m(K_Y+\Delta(f))\subset (\Omega_X^p)^{\otimes m}$, et donc que: 
$\psi^*(m(K_Y+\Delta(f))\subset mK_U$.

\

{\bf 7. LE ``COEUR".}

\

{\bf A. Le coeur.} 

\

{\bf 7.1. Th\'eor\`eme:} Soit $X$ (compacte K\" ahler connexe). Il 
existe une unique fibration $c_X:X\to C(X)$, d\'efinie \`a \'equivalence bim\'eromorphe pr\`es, telle que:

{\bf 1.} La fibre g\'en\'erale de $c_X$ est sp\'eciale.

{\bf 2.} $c_X$ est une fibration de type g\'en\'eral, ou constante (ssi $X$ est sp\'eciale, donc).

On appelle $c_X$ le {\bf coeur} de $X$, et $dim(C(X)):=ess(X)$ sa {\bf dimension essentielle}.

\

{\bf id\'ee de la d\'emonstration:} Unicit\'e de $c_X$ r\'esulte de {\bf 6.4}.

Existence. On proc\`ede par r\'ecurrence sur $n$. Si $X$ est sp\'eciale, $c_X$ est 
l'application constante. Sinon soit $f:X\to Y$ une fibration  de type g\'en\'eral avec $dim(Y)=p>0$ maximum. 
On veut montrer que $X_y$ est sp\'eciale. Soit $f=h\circ g$ le ``coeur relatif" de $f$, constitu\'e des deux fibrations $g:X\to Z$ et $h:Z\to Y$ telles que la restriction $g_y:=g_{\vert X_y}:X_y\to Z_y$ soit le coeur de $X_y$. (L'existence 
de ce coeur relatif est obtenue par la r\'ecurrence sur la dimension).

$$\xymatrix{X\ar[r]^g\ar[rd]_{f} & (Z)\ar[d]^h \\
& Y
}$$

\
Donc $g_y$ est de type g\'en\'eral ainsi que $h\circ g$. On d\'eduit de {\bf 5.5} que $g$ est aussi de type g\'en\'eral. 
Donc $dim(Z)=dim(Y)$, $Z=Y$, et $g=f$. Comme les fibres de $g$ sont sp\'eciales par construction, et que $g=f$, 
celles de $f$ le sont aussi. Donc $f=c_X$, par unicit\'e.

\

{\bf 7.2. Exemples:} 

\

{\bf 1.} $ess(X)=0$ ssi $X$ est sp\'eciale.

\

{\bf 2.} $ess(X)=n= dim(X)$ ssi $\kappa(X)=n$, ie: ssi $X$ est de type g\'en\'eral.

\

{\bf 3.} Si $X$ est une courbe, on a donc: $ess(X)=0$ ssi $g(X)=0,1$, et $ess(X)=1$ ssi $g(X)\geq 2$. On peut 
d\'ecrire: $c_X=(J_X\circ r_X)$ (voir {\bf 2.7) et {\bf 6.4}}.

{\bf 4.} On a fonctorialit\'e du coeur dans le sens suivant: si $f:X\to Y$ est une fibration, alors il 
existe une unique fibration $c_f:C(X)\to C(Y)$ telle que $c_f\circ c_X=c_Y\circ f$.

\

{\bf 5.} Si $f:X\to Y$ est une fibration dont la fibre g\'en\'erale $X_y$ est sp\'eciale, 
alors il existe $h:Y\to C(X)$ telle que 
$c_X=h\circ f$, par {\bf 6.4}. En particulier, si $f$ est soit $r_X:X\to R(X)$, le quotient rationnel de $X$, 
soit $J_X:X\to J(X)$, la fibration d'Iitaka-Moishezon de $X$ (d\'efinie si $\kappa(X)\geq 0$), 
on a une telle factorisation $h$. 
On dispose donc de factorisations $h_r:R(X)\to C(X)$, et $h_J:J(X)\to C(X)$ de 
$c_X=h_r\circ r_X=h_J\circ J_X$.

\

{\bf 6.} Lorsque $f=r_X$, il r\'esulte de [G-H-S] que $\Delta(r_X)=\emptyset$, et que $C(X)=C(R(X))$ (voir 
{\bf 4.3} et {\bf 2.11}).

\

{\bf 8.} On peut appliquer ce qui pr\'ec\`ede \`a $f:=J_{R(X)}:R(X)\to J(R(X))$ si $\kappa(R(X))\geq 0$ 
(la conjecture $-\infty$ affirme que c'est toujours le cas). 

\

On obtient ainsi une factorisation 
$h_J:J(R(X))\to C(X)=C(R(X))$ 
de $c_X:=h_J\circ (J_{R(X)}\circ r_X)$. 

\

$$\xymatrix{X\ar[r]^{r_X}\ar[rd]_{c_X} & R(X)\ar[r]^{J_{R(X)}}\ar[d]^{c_{R(X)}} & J(R(X))\ar[d]^{h_J} \\
& C(X)\ar[r]^{=} & C(R(X))
}$$

\

Cette factorisation est en fait le premier terme de la d\'ecomposition 
$c_X=(J\circ r)^n$ du coeur en une tour de fibrations \`a fibres soit avec $\kappa=0$ (au sens orbifolde), 
soit $\kappa_+=-\infty$ (une version faible de la connexit\'e rationnelle).

\

{\bf B. La d\'ecomposition du coeur.}

\

{\bf 7.5 D\'efinition:} $\kappa_+(X)=max\lbrace \kappa(Y,f), f:X\to Y\rbrace$, o\`u $f$ d\'ecrit 
l'ensemble des fibrations d\'efinies sur $X$. 

La d\'efinition lorsque $X$ est munie d'une structure 
d'orbifolde $(X/E)$ est similaire, en rempla\c cant 
$\kappa(Y,f)$ par $\kappa(Y,f,E)$ d\'efinie en {\bf 4.9}).

\

{\bf (7.6) Remarque:} Si $X$ est RC, $\kappa_+(X)=-\infty$, et cette condition 
implique que $X$ est $\kappa$-RC. 
Donc, si la conjecture {\bf (1.10)} est vraie, $X$ est RC ssi $\kappa_+(X)=-\infty$.

\

{\bf (7.7) Proposition:} Soit $f:(X/E)\to Y$ une fibration dont la base orbifolde et 
la fibre orbifolde g\'en\'erale ont $\kappa_+=-\infty$. Alors $\kappa_+(X/E)=-\infty$ aussi.

\

{\bf (7.8) Th\'eor\`eme:} Soit $(X/E)$ une orbifolde. Il existe une unique fibration 
$r=r_{(X/E)}:(X/E)\to Y$ telle que:

1. $\kappa(Y/\Delta(r,E))\geq 0$.

2. $\kappa_+(X_y/E_y)=-\infty$.

3. De plus: $p=dim(Y)<dim(X)=n$ ssi $\kappa(X)=-\infty$.

On appellera $r_{(X/E)}$ la ``$-\infty$-r\'eduction" de $(X/E)$.

\

{\bf 7.9.} Soit maintenant $f:X\to Y$ une fibration. On d\'efinit des fibrations $r$ et $J$ comme suit:

a. $r_f=r_{(Y/\Delta(f)}:Y\to R(f)$ dans tous les cas.

b. $J_f:=J_{(Y/\Delta(f))}:Y\to J(f)$ si $\kappa(Y/\Delta(f))\geq 0$ (voir {\bf (3.3)}).

\

$$\xymatrix{X\ar[r]^f & Y\ar[ld]_{r_f} \\
R(f)\ar[r]_{J_{r_f\circ f}} & J({r_f\circ f})
}$$

\

La compos\'ee $s_f:=J_{(r_f\circ f)}\circ r_f:Y\to S_f:=J({r_f\circ f})$ est donc
 bien d\'efinie dans tous les cas, par la condition {\bf 7.8.1}. Les fibrations $r_f,J_f$ coincident 
avec $r_X,J_X$ quand $f=id_X$.

De plus, $s_f=id_Y$ ssi: ou bien $Y$ est un point, ou bien $(Y/\Delta(f))$ est de type g\'en\'eral.

\

{\bf 7.10. Th\'eor\`eme:} Si $f:X\to Y$ est \`a fibre g\'en\'erale sp\'eciale, 
alors $s_f\circ f:X\to S_f$ aussi.

\

C'est une application facile d'une version orbifolde de {\bf (6.3)}. 

On en d\'eduit, 
par application it\'er\'ee des fibrations 
$s_f$, commen\c cant avec $f_0=id_X$ (qui est \` a fibres sp\'eciales!), une suite de fibrations 
$s^k:=(J\circ r)^k:X\to S_k(X), k\geq 0$, \`a 
fibres sp\'eciales. Cette suite stationne \`a $s^k:X\to S_k(X), k\leq n$, ssi: ou bien $S_k(X)$ 
est un point, ou bien $s^k$ une fibration de type g\'en\'eral.

On obtient donc:

\

{\bf 7.11. Th\'eor\`eme (de ``d\'ecomposition du coeur"):} Pour toute $X$, on a: $c_X=(J\circ r)^n$.

\

En particulier, si $X$ est sp\'eciale:

\

{\bf 7.12. Corollaire:} $X$ est sp\'eciale ssi $(J\circ r)^n$ est la fibration constante 
d\'efinie sur $X$ (ie: ssi $X$ peut \^etre d\'ecompos\'ee comme 
une tour de fibrations \`a fibres orbifoldes ayant soit $\kappa_+=-\infty$, soit $\kappa=0$).

\

{\bf 7.13. Remarque:} Les dimensions des $S_k(X), k\leq n$ d\'efinissent de nouveaux invariants bim\'ero-morphes 
de $X$, fonctoriels pour les applications m\'eromorphes surjectives. Ces invariants sont aussi, conjecturalement,
 des invariants de d\'eformation (voir {\bf 8.2}). 

On a un exemple simple de d\'ecomposition du coeur dans le cas des surfaces.

\

{\bf 7.14. Exemple:} Si $X$ est une surface. Alors: $c_X=(J\circ r)^2$, o\`u $J$ et $r$ 
sont des versions orbifoldes explicites
du quotient rationnel et de la fibration d'Iitaka-Moishezon. Plus pr\'ecis\'ement:

\

{\bf .} Si $\kappa(X)=-\infty$, alors $X$ est bim\'eromorphe \`a $\Bbb P^1\times C$, avec $g:=g(C)\geq 0$. Donc:

Si $g=0$, $r_X=c_X$ est l'application constante. 

Si $g=1$, $r_X:X\to C$ est la seconde projection. Et $J_C=J_{R(X)}=c_X$ est 
l'application constante. 

Si $g\geq 2$, $r_X:X\to C$ est la seconde projection, et $J_C=id_C$. Donc 
$c_X=J_{R(X)}\circ r_X:X\to C$. 

\

{\bf .} Si $\kappa(X)=0$, $J_X=J\circ r_X=c_X$ est l'application constante de $X$. 

\

{\bf .} Si $\kappa(X)=2$, on a: $c_X=r_X=J_X=id_X$, donc encore: $c_X=J_X\circ r_X$.

\

{\bf .} Si $\kappa(X)=1$, on a $r_X=id_X$, et $J:=J_X:X\to C$ est une fibration elliptique sur une courbe $C$. 
On fait intervenir la base orbifolde $B:=(C/\Delta(J))$ de $J$. On d\'efinit $r:=r_B$ et $J:=J_B$ en fonction 
de $\kappa(B)$ comme dans le cas non-orbifolde. (ie: $r_B$ est l'application constante ssi $\kappa(B)=-\infty$, 
et $J_B$ est l'application constante ssi $\kappa(B)=0$). On constate alors que 
$c_X=(J\circ r)^2$ dans tous les cas. 

C'est le cas le plus simple dans lequel les versions orbifoldes de $J$ et $r$ interviennent. 
On trouvera la decomposition (nettement plus compliqu\'ee) de $c_X$ en dimension $3$ dans [C01].

\

{\bf 8. CONJECTURES.}

\

{\bf A. D\'eformations.}

\

{\bf 8.1. Conjecture:} La classe des vari\'et\'es sp\'eciales est 
stable par d\'eformation et sp\'ecialisation.

(Ceci est vrai jusqu'en dimension $2$).

\

Plus g\'en\'eralement:

\

{\bf 8.2. Conjecture.} $ess(X)$ et $c_X$ sont invariants par 
d\'eformation de $X$. Plus g\'en\'eralement: 
toutes les fibrations interm\'ediaires 
$(J\circ r)^k$ et $r\circ (J\circ r)^k$ se d\'eforment avec $X$, pour $n\geq k\geq 0$. 

\

{\bf 8.3. Remarques:} L'invariance par d\'eformation de $\kappa$ 
(due dans des cas substantiels \`a Y.T.Siu) 
r\'esout affirmativement le cas $k=0$, lorsque $\kappa\geq 0$. Si la conjecture 
$-\infty$ est vraie, alors la conjecture {\bf 7.8} est vraie pour $r$ (ie: le cas $k=0$).

\

{\bf B. Le groupe fondamental.}

\

{\bf 8.4. Conjecture.} Si $X$ est sp\'eciale, alors $\pi_1(X)$ est presque-ab\'elien 
(ie: est ab\'elien, si l'on remplace $X$ par un rev\^etement \'etale fini ad\'equat). 

\ 

Cette conjecture est vraie si $X$ est RC (car $\pi_1(X)=1$, dans ce cas), ou 
si $\pi_1(X)$ est lin\'eaire 
(ie: plongeable dans un $Gl(n,\Bbb C)$), ou si $c_1(X)=0$, par Calabi-Yau. 
Un cas crucial ouvert est $\kappa(X)=0$.

Le cas g\'en\'eral se r\'eduit d'ailleurs aux versions orbifoldes 
des cas $\kappa_+=-\infty$ et $\kappa=0$, \`a 
l'aide de {\bf 7.7} ci-dessous. On peut en effet (voir [C01]) naturellement d\'efinir 
les notions de groupe fondamental d'une orbifolde (\`a multiplicit\'es ``standard" (voir {\bf 3.1})), et 
d'orbifolde sp\'eciale).

\

{\bf C. La pseudom\'etrique de Kobayashi.}

\

Si $X$ est une vari\'et\'e complexe compacte (K\" ahler, ici), 
on note $d_X:X\times X\to \Bbb R$ sa pseudom\'etrique de Kobayashi.

\

{\bf 8.5. Conjecture:} $X$ est sp\'eciale ssi $d_X\equiv 0$.

\

Cette conjecture est vraie si $X$ est RC, et pour $n=1$ et $n=2$ 
(sauf, peut-\^etre, si $\kappa(X)=2$, auquel cas elle est 
une version faible de la conjecture hyperbolique de Lang). Le cas crucial ouvert est, 
\`a nouveau: $\kappa=0$.

On peut d\'emontrer cette conjecture dans le cas particulier o\`u $X$ admet une fibration $f:X\to Y$ 
dont les fibres lisses sont RC ou des vari\'et\'es ab\'eliennes, pourvu que deux points 
g\'en\'eriques de la base puissent \^etre joints par 
l'image d'une courbe enti\`ere de $Y$ (ie: une application holomorphe de $\Bbb C$ dans $Y$).

\

Une autre situation dans laquelle on peut \'etablir que $d_X\equiv 0$ implique: 
$X$ est sp\'eciale, est celle du 
th\'eor\`eme {\bf 6.8}. En effet, $d_X$ s'annule alors sur l'adh\'erence de l'image de 
$\varphi:\Bbb C\to X$, suppos\'ee non d\'eg\'en\'er\'ee. Si 
cette image est dense, on a donc $d_X\equiv 0$, et $X$ est sp\'eciale par {\bf 6.8}. La situation de 
la remarque {\bf 6.9} est 
similaire, plus g\'en\'erale.

\

On peut d'ailleurs poser de nombreuses questions dans ce contexte. Parmi celles-ci: 

1. A-t'on $d_X\equiv 0$ si $d_X$ s'annule sur $U\times U$, o\`u $U$ est un ouvert non vide 
de $X$ (pour la topologie ``coh\'erente", bien s\^ur)?

2. Si $d_X\equiv 0$, existe-t'il une application holomorphe $\varphi:\Bbb C\to X$:

{\bf .} d'image dense dans $X$?

{\bf .} dont l'image contient un ensemble fini donn\'e, arbitraire, de $X$? (Une r\'eponse affirmative 
pour les couples d'\'el\'ements de $X$, avec variation alg\'ebrique des solutions avec les donn\'ees, 
entrainerait, par {\bf 6.9}, que $X$ est sp\'eciale si $d_X\equiv 0$).

\

Pour les $X$ arbitraires, le coeur permet de formuler la:

\

{\bf 8.6. Conjecture:} Soit $c_X:X\to C(X)$ le coeur de $X$, compacte K\" ahler. Alors:

{\bf 1.} Il existe une unique pseudom\'etrique $\delta_X$ sur $C(X)$ telle que $d_X=(c_X)^*(\delta_X)$.

{\bf 2.} $\delta_X=d_{C(X/\Delta(c_X))}$, la pseudom\'etrique de la base orbifolde de $c_X$.

{\bf 3.} $\delta_X$ est une m\'etrique sur un ouvert de Zariski dense $U$ de $C(X)$.

\

L'assertion {\bf 1.} n'est autre qu'une reformulation de la conjecture {\bf 8.5}. 

La condition {\bf 3.} est donc, compte tenu de {\bf 2}, la version orbifolde 
de la conjecture hyperbolique de Lang. 

La pseudom\'etrique de Kobayashi 
est d\'efinie sur une orbifolde $(Y/\Delta)$, si $\Delta=\sum_{j\in J}(1-1/m_j).D_j$, 
exactement comme dans le cas $\Delta=\emptyset$, mais 
en ne consid\'erant que les applications $h:\Bbb D\to Y$ du disque unit\'e $\Bbb D$ 
dans $Y$ telles que $h^{(k)}(z)$ soit tangente \`a tous les ordres $k>0$ non multiples de $m_j$ 
\`a $D_j$ en $h(z)\in D_j$ pour tout $z\in \Bbb D$ tel que $h(z)$ soit dans 
le lieu lisse du support de $\Delta$.

\

{\bf D. Arithm\'etique.} 

\

Soit $X$ une vari\'et\'e projective d\'efinie sur un corps de nombres $K$. Soit 
$K'\supset K$ un corps de nombres, et $X(K')\subset X$ l'ensemble des points $K'$-rationnels de $X$.

\

{\bf 8.7. Conjecture:}  $X$ est sp\'eciale ssi il existe un corps de  nombres 
$K'\supset K$ tel que l'ensemble $X(K')$ soit Zariski dense dans $X$. (On dit alors que $X$ est 
``potentiellement dense").

\

Cette conjecture n'est connue que dans des cas tr\`es particuliers (vari\'et\'es unirationnelles, vari\'et\'es 
ab\'eliennes), mais n'est pas m\^eme d\'emontr\'ee pour les surfaces $K3$ les plus g\'en\'erales, ou les vari\'et\'es RC 
(la version ``corps de fonctions" est cependant connue pour les vari\'et\'es RC, par [G-H-S]).

\

Remarquons que cette conjecture diff\`ere de celle de Colliot-Th\'el\`ene-Harris, qui affirme que 
$X$ est potentiellement dense ssi aucun rev\^etement \'etale 
fini $X'$ de $X$  n'admet de fibration $f:X\to Y$ avec $Y$ de type g\'en\'eral et 
$dim(Y)>0$. [B-T] construit 
en effet des exemples de $X$ simplement connexes n'admettant 
pas de telle fibration, et non sp\'eciales. 
(Autrement dit: les fibres multiples du coeur ne 
peuvent \^etre \'elimin\'ees par un rev\^etement \'etale de $X$ 
dans ce cas).

\

 En analogie avec {\bf 8.6}:

\

{\bf 8.8. Conjecture.} Soit $X$ une vari\'et\'e 
projective complexe d\'efinie sur un corps de nombres $K$. Alors $c_X(X(K))\cap U$ est fini, si 
$U$ est l'ouvert de Zariski dense de $X$ introduit en {\bf 8.6.3}. 

\

Autrement dit: $X(K)$ est concentr\'e sur un nombre fini de fibres de $c_X$ au-dessus de $U$.

\

Cette conjecture n'est connue que dans des cas tr\`es particuliers: 
courbes et sous-vari\'et\'es de vari\'et\'es ab\'eliennes.

\

{\bf 8.9 Remarque:} On pourra trouver la version ``corps de fonctions" 
des conjectures pr\'ec\'edentes dans [C01].

\

\

{\bf Bibliographie:}

[B-T]F.Bogomolov-Y.Tschinkel. Special Elliptic Fibrations. Preprint (2003).

[C81]F.Campana.Cor\'eduction alg\'ebrique d'un espace
analytique faiblement K\"ahl\'erien compact.Inv.\ Math.\ (1981),187-223.

[C91] F.Campana.Twistor spaces of class $\cal C$.J.Diff.Geom. 33 (1991),541-549. 

[C92]  F.~Campana.Connexit\'e rationnelle des
vari\'et\'es de Fano. Ann.\ Sc.\ ENS.\ 25 (1992), 539-45.

[C01] F.~Campana.Special Varieties and Classification Theory.
Math. AG/0110051. (A paraitre aux Ann. Inst. Fourier (2004)).

[De97]O.Debarre. Vari\'et\'es de Fano. S\'eminaire Bourbaki 1996/97, n¡827.

[De02]O.Debarre. Vari\'et\'es Rationnellement connexes. S\'eminaire Bourbaki 2001/02, n¡905.

[G-H-S] T.~Graber-M.~Harris-J.~Starr.Families of rationally
connected varieties. Preprint 2001.

[I] S.~Iitaka.Genera and Classification of Algebraic
Varieties. Sugaku 24 (1972), 14-27.

[K] Y.~Kawamata. Characterisation of Abelian
varieties.Comp.\ Math.\ (1981), 253-76.

[K-O] S.~Kobayashi-T.~Ochiai.Meromorphic mappings into
compact complex spaces of general type. Inv.\ Math.\ 31 (1975), 7-16.

[Ko-Mi-Mo]J.~Koll\'ar-Y.~Miyaoka-S.~Mori.Rationally
connected Varieties.J.\ Alg.\ Geom.\ 1 (1992), 429-448.

[L]S.Lang.Hyperbolic and Diophantine
Analysis. Bull.\ AMS 14(1986), 159-205.

[Mi]Y.Miyaoka.On the Kodaira Dimension of a Minimal Threefold. Math. Ann. 281 (1988),325-332.

[M]S.Mori.Flip Theorem and the existence of Minimal odels for Threefolds. J.AMS. 1 (1988), 117-253.

[Mo]B.~Moishezon.Algebraic Varieties and Compact
Complex Spaces.Actes du Congr\` es International des Math\'ematiciens, Nice 1970.
Vol.~2, 643-648.

[Sh]I.Shafarevitch. Algebraic Surfaces.Proceedings of the Stekhlov Institute. vol 75. 1965.

[U] K.~Ueno.Classification Theory of Algebraic Varieties and
Compact Complex Manifolds.LNM 439 (1975), Springer Verlag.

[V] E.~Viehweg.Vanishing theorems and and positivity
of Algebraic fibre spaces.J.\ Proc.\ Int.\ Congr.\ Math.\
Berkeley (1986).

\end

$$\xymatrix{
&M\ar@{-}[ldd]_3\ar@{-}[rdd]^2&\\
&L\ar@{--}[ld]\ar@{--}[rd]&\\
\Bbb{Q}[b]\ar@{-}[rd]_2&&\Bbb{Q}[c]\ar@{-}[ld]^3\\
&\Bbb{Q}&}$$

\

\